\documentclass{conm-p-l} 
\usepackage{amssymb}

\newcommand{\R}{{\mathbb R}}
\newcommand{\ep}{\epsilon}
\newcommand{\pd}{\partial}
\newcommand{\comp}{\operatorname{comp}}

\theoremstyle{plain}
\newtheorem{theorem}{Theorem}[section]

\newtheorem{lemma}[theorem]{Lemma}

\theoremstyle{definition}
\newtheorem{example}[theorem]{Example}

\theoremstyle{remark}
\newtheorem{remark}[theorem]{Remark}

\numberwithin{equation}{section}

\begin{document}

\title[Melnikov Method]{Melnikov Method for Autonomous Hamiltonians}
\author{Clark Robinson}
\address{Department of Mathematics, Northwestern University, Evanston IL
60201}
\email{clark@math.nwu.edu}
\thanks{This research partially supported by grants from the 
National Science Foundation.}

\begin{abstract}
This paper presents the method of applying the Melnikov method to
autonomous Hamiltonian systems in dimension four.  Besides giving an
application to Celestial Mechanics, it discusses the problem of convergence
of the Melnikov function and the derivative of the Melnikov function.
\end{abstract}

\maketitle

\section {Introduction}  
The original Melnikov method applies to periodically forced Hamiltonian
systems in dimension two (one degree of freedom).  
In this paper, our main emphasis is to show how a modification of this
method applies directly to  autonomous Hamiltonian systems in dimension four 
(two degrees of freedom).  
These original unperturbed system needs to be completely integrable.
By restricting to the energy surface,
it is possible to convert the autonomous system in dimension four into a 
periodically forced Hamiltonian
systems in dimension two and then apply the usual Melnikov method.
It is our contention that this is not the natural method to treat this type
of system.  In fact the calculations are much more direct and natural 
when done directly with the autonomous system.  
Thus in this paper we advocate applying the Melnikov method
directly to the autonomous Hamiltonian system.
These ideas go back at least to the papers of P. Holmes and J. Marsden,
\cite{hm82} and \cite{hm83}. Also see \cite{CR-mel} and \cite{lu}.
In Section 3 we discuss the main theorems in the theory of autonomous
systems and their proofs.  In Section 4 we indicate how these theorems can
be applied to the circular restricted three body problem.

Because we want to connect our treatment with the material for periodically
forced systems, we start by describing the periodically forced
situation using our approach in Section 2.

This paper also points out the connection between the convergence
of the improper integral giving the Melnikov function and the smoothness
and dependence on parameters of the stable manifold of a periodic orbit.  
(See the discussion of Lemma \ref{lemma} and the remarks which follow it.)
This connection has also been made recently in the context of Celestial
Mechanics by J. Casasayas, E. Fontich, and A. Nunes, \cite{cfn94}.

A final point made in this paper is the possibility of differentiating the 
Melnikov function by differentiating the integrand of its representation as
a  improper integral.  
See Theorem \ref{thm-deriv} in Section 3.

\section{Periodically force systems}  

Consider a Hamiltonian $H:\R^2 \to \R$, with associated Hamiltonian vector
field $X = X_H$. Next, consider the periodically force perturbation,
\begin{equation*}
\begin{split}
\dot \mathbf{z} &= X(\mathbf{z}) + \ep Y(\mathbf{z},\tau)\\
\dot \tau       &= 1
\end{split} 
\end{equation*}
where $\tau$ is considered as a periodic variable of some period $T$.
We assume that $X$ has a homoclinic orbit to a hyperbolic fixed point (in the
$\mathbf{z}$-space).  Let $\gamma_0$ be the associated periodic orbit in
the $(\mathbf{z},\tau)$-space, $\tilde S$, and 
\begin{equation*}
\Gamma_0 = \{(\mathbf{z}_0,\tau_0) \:
\mathbf{z}_0 \text{ is homoclinic to } \gamma \}
\end{equation*}
 
be the set of homoclinic orbits in the extended phase space $\tilde S$.

Since $\gamma_0$ is hyperbolic, it can be continued for $\ep \neq 0$
to a periodic orbit $\gamma_\ep$ in the extended phase space $\tilde S$.
For $\sigma$ equal $s$ and $u$,
let $W^\sigma(\gamma_\ep, X+\ep Y)$ be the stable and unstable manifolds of
$\gamma_\ep$ respectively.
In order to determine the orbits on $\Gamma_0$ to a particular
orbit on $W^\sigma(\gamma_\ep, X+\ep Y)$  for $\ep \neq 0$, 
we let $N_{(\mathbf{z}_0,\tau_0)}$ be a (one dimensional) transversal to 
$\Gamma_0$ at a point $(\mathbf{z}_0,\tau_0) \in \Gamma_0$.  
Then for $\sigma = s,u$, set
\begin{equation*}
\zeta^\sigma(\mathbf{z}_0,\tau_0,\ep) = W^\sigma(\gamma_\ep, X+\ep Y)
\cap N_{(\mathbf{z}_0,\tau_0)}.
\end{equation*}

In our development, we use strongly the fact that these stable and unstable
manifold depend differentiably on the parameter and position on the
manifold, so are smooth in the space $\tilde S$, i.e., for 
$\sigma = s$ and $u$,  
$\zeta^\sigma(\mathbf{z}_0,\tau_0,\ep)$ is a smooth function of all its
variables.
 
The Hamiltonian function is constant on $\Gamma_0$ and has nonzero
gradient, so it is a good measure of the displacement
perpendicular to $\Gamma_0$. 
With this motivation, we define the \emph{Melnikov function} by
\begin{equation*}
M(\mathbf{z}_0,\tau_0) 
=\frac{\pd}{\pd \ep} \big[H \circ \zeta^u(\mathbf{z}_0,\tau_0,\ep)
- H \circ \zeta^s(\mathbf{z}_0,\tau_0,\ep)\big]\big|_{\ep =0}.
\end{equation*}
Thus, $M$ is a measurement of the infinitesimal separation of the stable
and unstable manifolds as measured by the Hamiltonian function $H$.
By parameterizing the solutions of the equations and restricting to a
transverse cross section, $M$ can be considered a function of a real
variable (as it usually is).  
Notice that we use the fact that the stable and unstable manifolds depend
smoothly on the parameter to define $M$.

There are two theorems related to $M$.  The first states that a
nondegenerate zero of $M$ corresponds to a transverse homoclinic orbit of
$\gamma_\ep$ for $\ep \neq 0$ but small enough.  We refer to
\cite{CR-mel}, \cite{wiggins}, or other standard 
references for this theorem.  
(Also see Theorem \ref{thm-trans} in the next section.)
The second theorem relates $M$ to an improper integral. 
We state this second
theorem in order to emphasis the fact that with the proper statement, it
can be calculated in any coordinate system (not just symplectic
coordinates).

\begin{theorem}\label{thm-mel-per}
\begin{align}
M(\mathbf{z}_0,\tau_0) 
&= \int_{-\infty}^\infty \mathbf{k}\cdot (X \times
Y)_{(\phi_0^t(\mathbf{z}_0), t+\tau_0)}\, dt \tag{*}\label{eq:cross} \\
&= \int_{-\infty}^\infty (DH \cdot Y)_{(\phi_0^t(\mathbf{z}_0),
t+\tau_0)}\, dt. \tag{**}\label{eq:ham}
\end{align}
\end{theorem}

We remark that the integral in equality (\ref{eq:cross}) is the usual 
one derived for the Melnikov function: 
the integrand is the scalar cross product (component in the 
$\mathbf{k}$ direction) or wedge product of the two vector fields.  
This integrand must be calculated in symplectic coordinates.  
The the integral in equality (\ref{eq:ham}) is
valid in any set of coordinates, even those which are not symplectic.
See \cite{CR-mel} or \cite{wiggins} for a proof of the second form 
of the integral.  
Also in the discussion below of the autonomous case, we indicate a sketch
of the proof why the above theorem is true.
The proof in \cite{CR-mel} uses the smoothness of the stable and unstable
manifolds and not the uniform estimate for $0 \leq t < \infty$ 
on the difference of solutions for $\ep \neq 0$ and $\ep =0$
as many proofs do. (See e.g. \cite{cfn93} and \cite{cfn94}.)

\section{Autonomous Hamiltonians in four dimensions}

In this section, the phase space is some four dimensional space which we
label as $S^4$: in Celestial Mechanics applications $S^4$ is usually $\R^4$ or 
some product of copies of the circle (periodic variables) and copies of the
reals.  The autonomous Hamiltonian $H^\ep:S^4 \to \R$ is given as
\begin{equation*}
H^\ep(\mathbf{z}) = H_0(\mathbf{z}) + \ep H_1(\mathbf{z}) +O(\ep^2).
\end{equation*}
The corresponding Hamiltonian vector field is given as
\begin{equation*}
X^\ep(\mathbf{z}) = X_0(\mathbf{z}) + \ep Y(\mathbf{z}) +O(\ep^2).
\end{equation*}

The basic assumptions are as follows.
\begin{enumerate}
\item[(A1)] We assume that $X_0$ is completely integrable, with a second
independent integral $F$.  (The function $F$ is mostly likely not an
integral of $X^\ep$ for $\ep \neq 0$.)
\item[(A2)] We assume that $X_0$ has a hyperbolic periodic orbit $\gamma_0$ 
with set of homoclinic  orbits 
\begin{align*}
\Gamma_0 &= \comp[W^u(\gamma_0,X_0) \setminus \gamma_0] \\
&= \comp[W^s(\gamma_0,X_0) \setminus \gamma_0].
\end{align*}
Here $\Gamma_0$ is one of the components of the stable and unstable manifolds 
indicated and has dimension two because of the second integral.
\end{enumerate}

\begin{example}  
The easiest example of the above type of system is one which 
decouples for $\ep =0$ into a subsystem with a homoclinic orbit and another 
That is oscillator: e.g.,   assume 
$H_0(\mathbf{z},\mathbf{w}) = F(\mathbf{z}) + G(\mathbf{w})$
where
\begin{equation*}
F(z_1,z_2) = \frac{z_2^2}{2} - \frac{z_1^2}{2} + \frac{z_1^3}{3}
\end{equation*}
has a homoclinic orbit to the fixed point at $\mathbf{0}$, and
\begin{equation*}
G(w_1,w_2) = \frac{\alpha}{2}(w_1^2 + w_2^2)
\end{equation*}
is an oscillator.  A simple coupling is of the type 
$H_1(\mathbf{z},\mathbf{w}) = z_1 w_1$.
\end{example}

\begin{example}  
P. Holmes and J. Marsden \cite{hm83}  gave an example of a 
system which does not decouple. 
(Compare this example with the closely related one from Celestial 
Mechanics given in \cite{cfn93}.)  
Let
\begin{equation*}
H^\ep(p,q,\theta, I) = \frac{p^2}{2} + U(q) + \frac{I^2}{2 q^2} 
+ \ep \sin(\theta),
\end{equation*}
so 
\begin{equation*}
\dot I = -\ep \cos(\theta)
\end{equation*}
and $F = I$ is a second integral for $\ep = 0$.
The potential function $U(q)$ must be chosen so that there is a homoclinic 
orbit in the $(p,q)$-subsystem for $\ep = 0$ and $I$ fixed.
\end{example}

We set up the definition of the Melnikov function using transversals 
much as we did above for the periodically forced case.
For each point $\mathbf{z_0} \in \Gamma_0$, let $N_{\mathbf{z_0}}$
be a two dimensional plane perpendicular to $\Gamma_0$ at $\mathbf{z_0}$.
Let $\gamma_\ep$ be the continuation of the periodic orbit $\gamma_0$ for 
$\ep \neq 0$ with $\gamma_\ep \subset (H^\ep)^{-1}(h_0)$
where $h_0 = H_0(\gamma_0)$ is the energy of the unperturbed orbit.  
Then for $\sigma = s$ and $u$ and $\ep \neq 0$, let
\begin{equation*}
\zeta^\sigma(\mathbf{z_0},h_0,\ep) = N_{\mathbf{z_0}} \cap
W^\sigma(\gamma_\ep, X^\ep)
\end{equation*}
be the continuation of $\mathbf{z_0}$.
Notice that both the stable and the unstable manifolds continue to lie in the 
same energy surface $(H^\ep)^{-1}(h_0)$, 
so $H_0$ is not a good measure of
the separation;  on the other hand, the second integral $F$ is constant 
on $\Gamma_0$ and
is independent of $H_0$ so it is a good measure of the distance 
perpendicular to $\Gamma_0$ within this energy surface.
The {\emph{Melnikov function}}, $M:\Gamma_0 \to \R$ is then 
defined as the infinitesimal
displacement of the stable and unstable manifolds as measured by the second
integral $F$:
\begin{equation*}
M(\mathbf{z_0},h_0) 
= \frac{\pd}{\pd \ep} \big[F\circ\zeta^u(\mathbf{z_0},h_0,\ep)
- F\circ\zeta^s(\mathbf{z_0},h_0,\ep)\big]\big|_{\ep=0}.
\end{equation*}
By expressing $M$ in terms of specific coordinates and restricting to a 
transverse section, it is possible to give it
as a function of a single variable (plus the parameter $h_0$).
See the example in the next section.

Again, there are two theorems about the Melnikov function.  The first states
that a nondegenerate zero of $M$ corresponds to a transverse homoclinic point.
The second theorem indicates how to calculate $M$ in terms of an improper
integral.
In general for this case, the integral representation of $M$ is 
only conditionally convergent, so the times must be chosen correctly 
in the improper integral.  
The third theorem below justifies differentiating under the
integral sign to calculate the derivative of the Melnikov function.

\begin{theorem}\label{thm-trans}
Make the basic assumptions A1-A2).
If $M(\cdot, h_0):\Gamma_0 \to \R$ has a nondegenerate zero at 
$\mathbf{z_0} \in 
\Gamma_0$, then $X^\ep$ has a transverse homoclinic orbit 
$\mathbf{z_\ep}$ nearby for small enough nonzero $\ep$.
\end{theorem}

\begin{theorem}\label{thm-mel}
Make the basic assumptions (A1-A2).
Assume that $T_j$ and $-T_j^*$ are chosen so that
the points
$\phi_0(T_j,\mathbf{z_0})$ and $\phi_0(-T_j^*,\mathbf{z_0})$
converge to the same point $\mathbf{z_\infty} \in \gamma_0$. Then
\begin{equation*}
M(\mathbf{z_0},h_0) 
= \lim_{j \to \infty} \int_{-T_j^*}^{T_j}\,
 (DF \cdot Y)_{\phi_0(t,\mathbf{z_0})}\, dt. \label{eq:thm-mel}
\end{equation*}
\end{theorem}

\begin{theorem}\label{thm-deriv}
Make the basic assumptions (A1-A2) as in Theorem \ref{thm-mel} 
but assume that the integral converges absolutely (not conditional 
convergence). 
(See remarks following Lemma \ref{lemma} for conditions which imply this
assumption.)
Then
\begin{equation*}
D_{\mathbf{z}}M_{(\mathbf{z},h_0)} \mathbf{v}
= \int_{-\infty}^{\infty}\, D_{\mathbf{z}}\big( 
(DF \cdot Y)_{\phi_0(t,\mathbf{z})}\big)\mathbf{v} \, dt. 
\label{eq:thm-deriv}
\end{equation*}
In the above equation, $D_z$ is differentiation with respect to the initial
spatial coordinate $\mathbf{z}$, and is calculated in the direction of the
vector $\mathbf{v}$.
\end{theorem}

\begin{remark}  The integral in Theorems \ref{thm-mel} and \ref{thm-deriv}
can be calculated in any 
coordinates, even non-symplectic coordinates.  This is relevant to applications
in Celestial Mechanics because McGehee coordinates are non-symplectic 
coordinates.   
See Section 4.
\end{remark}

\begin{remark}  Although the convergence is only conditional in general, there
are several circumstances in which the convergence is absolute.  We make this 
more explicit below.  
\end{remark}

The proof of Theorem \ref{thm-trans} follows as usual from the 
Implicit Function Theorem, so we do not give any details.
See \cite{wiggins}.

\begin{proof}[Idea of the proof of Theorem \ref{thm-mel}]
As in the usual proof, 
\begin{equation*}
M(\mathbf{z_0},h_0) = \int_{-T_j^*}^{T_j}\,
 (DF \cdot Y)_{\phi_0(t,\mathbf{z_0})}\, dt
+ R(-T_j^*, T_j)
\end{equation*}
where the remainder term 
\begin{equation*}
\begin{split}
R(-T_j^*, T_j) =
&\frac{\pd}{\pd\ep}\big[F\circ\phi(-T_j^*,\zeta^u(\mathbf{z_0},\ep),\ep)\big]
\big|_{\ep=0} \\
&- \frac{\pd}{\pd\ep}\big[F\circ\phi(T_j,\zeta^s(\mathbf{z_0},\ep),\ep)\big]
\big|_{\ep=0}.
\end{split}
\end{equation*}
In \cite{CR-mel}, we prove a lemma about the convergence of the 
remainder. 
The statement of Lemma 2.1 in \cite{CR-mel} only includes the final
statement of the following lemma but the proof in that paper implies the
following extended statement.
Clearly Lemma \ref{lemma} implies the Theorem \ref{thm-mel}.  
\end{proof}

\begin{lemma}\label{lemma} 
Assume that $T_j$ and $-T_j^*$ are chosen so that
$\phi_0(-T_j^*,\mathbf{z_0})$ converges to $\mathbf{z^-_\infty}$
and $\phi_0(T_j,\mathbf{z_0})$ converges to $\mathbf{z^+_\infty}$
as $j \to \infty$. 
Then 
\begin{equation*}
\frac{\pd}{\pd\ep}\big[F\circ\phi(-T_j^*,\zeta^u(\mathbf{z_0},\ep),\ep)\big]
\big|_{\ep=0}
\end{equation*}
converges to 
\begin{equation*}
\frac{\pd}{\pd\ep}\big[F\circ\gamma(\mathbf{z^-_\infty},\ep)\big]
\big|_{\ep=0}
\end{equation*}
and 
\begin{equation*}
\frac{\pd}{\pd\ep}\big[F\circ\phi(T_j,\zeta^s(\mathbf{z_0},\ep),\ep)\big]
\big|_{\ep=0}
\end{equation*}
converges to 
\begin{equation*}
\frac{\pd}{\pd\ep}\big[F\circ\gamma(\mathbf{z^+_\infty},\ep)\big]
\big|_{\ep=0}.
\end{equation*}
Here $\gamma(\mathbf{z},\ep)$ is the function which gives the perturbed
orbit $\gamma_\ep$ for $\ep \neq 0$ as a graph over the unperturbed orbit.
In particular, if $\mathbf{z^-_\infty} = \mathbf{z^+_\infty}$
then $R(-T_j^*,T_j)$ goes to zero as $j$ goes to infinity.
\end{lemma}

\begin{remark}
By Lemma \ref{lemma},
the improper integral in Theorem \ref{eq:thm-mel}  always gives
$M(\mathbf{z_0},h_0)$ provided $T_j$ and $-T_j^*$ are chosen so that 
$\mathbf{z^-_\infty} = \mathbf{z^+_\infty}$.
\end{remark}

\begin{remark}
If $DF_{(\mathbf{z})} \equiv 0$ for all $\mathbf{z} \in \gamma_0$,
then the convergence is not conditional and the integral is just the
improper integral from $-\infty$ to $\infty$.
This feature is essentially what is true in the
proof for the usual Melnikov function when using the function $H$ which 
vanishes at the fixed point.
\end{remark}

\begin{remark}
If the periodic orbit does not move with $\ep$,
$\gamma_\ep \equiv \gamma_0$,
then the convergence is also a usual improper integral (and not conditional)
since the limiting values are zero and do not need to cancel each other.   
When applying this theorem to parabolic orbits which are asymptotic to
periodic orbits at infinity, 
this is often the case since the orbits at 
infinity are not affected by the small coupling terms.
See the example of the circular restricted three body problem in Section 4.
\end{remark}

\begin{remark}
What makes the convergence work in the proof of Lemma \ref{lemma} is 
the stable manifold theory.
Casa\-sayas, Fontich, and Nunes proved in \cite{cfn94} that this type of
result is often true even when the periodic orbit is weakly hyperbolic
(from higher order terms).  
Notice that they emphasize that it is possible to prove a stable manifold 
theorem in the parabolic case with a particular smoothness. 
Then then prove that this smoothness of the stable manifold is sufficient
to imply a convergence of 
the above type to give the Melnikov function as an improper integral.
They treat the usual time periodic perturbations
but their results go over to the above context of autonomous systems
as well.  
\end{remark}

\begin{proof}[Proof of Theorem \ref{thm-deriv}]
The fact that we can differentiate under the integral sign follows from
a theorem in analysis about differentiating improper integrals.
What is needed to differentiate under the integral sign of an improper
integral is that the integrals for $M(\mathbf{z},h_0)$ and 
$D_{\mathbf{z}}M_{(\mathbf{z},h_0)} \mathbf v$ converge uniformly in
$\mathbf{z}$. 
See \cite{lang} or \cite{lewin}.
In our situation the convergence of both integrals is uniform for 
the same reason
that the integral converges in the proof of 
Lemma \ref{lemma} as given in  \cite{CR-mel}.  
\end{proof}

\section{Planar Restricted Three Body Problem Example}

Both J. Llibre and C. Simo in \cite{ls} and Z. Xia in \cite{xia-mel} 
applied the usual Melnikov method to the circular planar three body problem.
In this section, we indicate how the same results can be obtain (possibly)
more directly by using the method for the autonomous systems indicated in
the last section.  Our remarks about the applicability in non-symplectic 
coordinates also justifies making the calculations in McGehee coordinates.

Below we sketch the introduction of the coordinates into the problem.
See the above references for more detailed explanations and treatment.
Our treatment most closely follows that in \cite{xia-mel}.

In this situation, the third mass $m_3 \equiv 0$, $m_1 = \mu$ is a small
parameter, and $m_2 = 1 - \mu$.  
Jacobi coordinates are used as a starting place.
The problem is circular so the variable
\begin{equation*}
\mathbf{Q} = (\cos(t), \sin(t))
\end{equation*}
gives the motion of the vector from $m_2$ to $m_1$.
The position of $m_3$ relative to $m_1$ and $m_2$ is given by $\mathbf{q}$.

Next, McGehee coordinates are introduced.  The variable $x$ measures the 
distance from infinity, 
\begin{equation*}
x^{-2} = |\mathbf{q}|,
\end{equation*}
so $x$ goes to zero as $|\mathbf{q}|$ goes to infinity.
The direction to $m_3$ is given by $\theta$,
\begin{equation*}
 \frac{\mathbf{q}}{|\mathbf{q}|} = (\cos \theta, \sin \theta ).
\end{equation*}
This angle is measured in rotating coordinates by introducing the 
$2\pi$-periodic variable $s = t - \theta$.
The momentum variable for $m_3$, $\mathbf{p}$ is broken down into radial and 
angular components by introducing the variables $y$ and $\rho$ respectively
by the equation
\begin{equation*}
\mathbf{p} = [y + x^2 \rho i] [\cos \theta + i \sin \theta ].
\end{equation*}
(This last equations uses the notation of complex variables for points in 
$\R^2$.)
The energy equation in these variables becomes
\begin{equation*}
H^\mu = \frac{1}{2} y^2 + \frac{1}{4} x^4 \rho^2 - U^\mu
\end{equation*}
where $U^0 = x^2$. 

When $\mu = 0$, $\dot \rho = 0$, so $\rho$, which is the angular momentum
of $m_3$ (the Jacobi constant), is a second integral.
After a change of the time variable, the equations become
\begin{align*}
x' &= -\frac{1}{2} x^3 y \\
y' &= -x^4 + x^2\rho^2 + \mu g_1(x,s) + O(\mu^2) \\
\rho' &= \mu g_2(x,s) + O(\mu^2) \\
s' &= 1 - x^4 \rho.
\end{align*}
Here $s$ is considered a $2\pi$ periodic variable.

When $\mu = 0$, $\begin{pmatrix} x \\ y \end{pmatrix} = 
\begin{pmatrix} 0 \\ 0 \end{pmatrix}$
is a fixed point in $\begin{pmatrix} x \\ y \end{pmatrix}$-space.
This fixed point corresponds to a periodic orbit in the total space
(with $s$ the periodic variable).
This periodic orbit is weakly hyperbolic but it has smooth
stable and unstable manifolds which correspond to the parabolic orbits of the
original problem.  We parameterize the homoclinic orbits 
\begin{align*}
x &= \xi_0(t,\rho_0,s_0) \\
y &= \eta(t,\rho_0,s_0)  \\
s &= s_0 + \sigma(t,\rho_0)
\end{align*}
where we take the initial conditions so that $\eta(0,\rho_0,s_0) = 0$,
i.e. the orbits cross the transverse section $y = 0$ at $t = 0$.  
Fixing $\rho_0$ at $t = 0$ is like fixing the energy of the system.

For $\mu = 0$, $\rho$ is constant and so is a second integral of the
system. 
Thus the Melnikov function is
\begin{align*}
M(s_0,\rho_0)
&= \int_{-\infty}^{\infty}\ (D\rho)\cdot Y_{\phi_0(t,\rho_0,s_0)}\, dt \\
&= \int_{-\infty}^{\infty}\ g_2 \circ \phi_0(t,\rho_0,s_0)\, dt  \\
&= \int_{-\infty}^{\infty}\ \xi_0^4(t)\sin(s_0+\sigma(t))
S(t,\rho_0,s_0,) \ dt,
\end{align*} 
where
\begin{equation*}
S(t,\rho_0,s_0,) = \big[1 - \frac{1}{[1+2\xi_0^2(t)\cos(s_0+\sigma(t))
+\xi_0^4(t)]^{3/2}\big]}.
\end{equation*}
The above form of the Melnikov function is derived in \cite{xia-mel}
by starting with the form arising from  a periodically forced system  
and then integrating by parts.  
This type of derivation is not surprising since the equivalence of tbe 
usual Melnikov integral and that obtained for autonomous sytems
is proved in \cite{CR-mel} using integration by parts. 
The point that is being made here is the above treatment is very direct 
and seems less ad hoc then what is found in \cite{xia-mel}.
(I am certain that Z. Xia was very aware of the representation of the 
Melnikov integral by means of a second integral and that he just 
preferred using the more familiar form.)
Notice that the orbit at infinity does not vary with $\mu$ so the 
integral converges absolutely.

Once the above form is obtained, symmetry considerations imply that there
is a zero of $M$ for $s_0 = 0$ and $\pi$, 
$M(0,\rho_0) = M(\pi,\rho_0) = 0$. 

The orbit at infinity does not vary with $\mu$, so the above integral
representation of $M$ is not conditional.  Therefore by Theorem 
\ref{thm-deriv}, $\dfrac{\pd M}{\pd s_0}(s_0, \rho_0)$ can be calculated by
differentiating under the integral sign with respect to $s_0$ which
appears in a very simple way. (The integral of this derivative with respect
to $s_0$ also converges uniformly in $s_0$.)
Special arguments are needed to get that the improper integral
representing $\dfrac{\pd M}{\pd s_0}(s_0, \rho_0)$
for $s_0 =0$ or $\pi$ is nonzero, and so the corresponding zero of 
$M(\cdot,\rho_0)$ is a nondegenerate zero as a function of $s_0$.
These arguments which verify this fact done very differently 
in \cite{xia-mel} and \cite{ls}.  
 
\bibliographystyle{amsplain}
\bibliography{mel}
\end{document}